\documentclass[twocolumn]{autart}    
\usepackage{amssymb,amsmath}
\usepackage{graphicx}  
\usepackage[usenames,dvipsnames]{xcolor}
\usepackage{bbold}
\usepackage{units}

\usepackage{pgfplots}
\pgfplotsset{compat=newest}
\usepackage{pgfmath,pgffor,pgfplotstable}
\usetikzlibrary{calc,decorations.markings,intersections}
\usepgfplotslibrary{fillbetween,patchplots,colorbrewer}

\usepackage{tikzscale}
\usetikzlibrary{external}
\newtheorem{theorem}{Theorem}

\newtheorem{lemma}{Lemma}
\newtheorem{remark}{Remark}
\newtheorem{assumption}{Assumption}
\newtheorem{proposition}{Proposition}
\newtheorem{definition}{Definition}
\newtheorem{example}{Example}

\newcommand{\tr}{\intercal}
\newcommand{\defeq}{:=}
\newcommand{\f}{{\sf{f}}}
\newcommand{\ep}{{\sf{e}}}
\newcommand{\vp}{{\sf{v}}}
\newcommand{\one}{\mathbf{1}}

\usepackage{xpatch}
\makeatletter 
\xpatchcmd{\runningauthor@fmt}{\global\edef}{\protected@xdef}{}{}
\xpatchcmd{\runningauthor@fmt}{\global\edef}{\protected@xdef}{}{}
\xpatchcmd{\author@fmt}{\edef}{\protected@edef}{}{}
\def\@xnamedef#1{\expandafter\protected@xdef\csname #1\endcsname}
\def\ead@au#1{\protected@edef\@ead@au{#1}}
\def\add@xtok#1#2{\begingroup
  \protected@xdef\@act{\global\noexpand#1{\the#1#2}}\@act
\endgroup}
\def\no@harm{}
\makeatother

\newcommand{\legendbox}[1]{%
  \textcolor[RGB]#1{\rule{\fontcharht\font`X}{\fontcharht\font`X}}%
}

\begin{document}

\begin{frontmatter}

\title{Control Lyapunov Function Design via Configuration-Constrained Polyhedral Computing}

\thanks[footnoteinfo]{This paper was not presented at any IFAC 
meeting. Corresponding author M.~E. Villanueva.}

\author[ST]{Boris Houska}
\ead{borish@shanghaitech.edu.cn}, 
\author[LUH]{Matthias A.~M{\"u}ller}
\ead{mueller@irt.uni-hannover.de},
\author[IMT]{Mario E.~Villanueva}
\ead{me.villanueva@imtlucca.it}

\address[ST]{ShanghaiTech University, China}%
\address[LUH]{Leibniz University Hannover, Germany}%
\address[IMT]{IMT School for Advanced Studies Lucca, Italy}%

\begin{keyword}
Linear Systems, Control Lyapunov Functions, Polyhedral Computing, Convex Optimization, Model Predictive Control, Uncertain Control Systems, Min-max 
Model Predictive Control
\end{keyword}

\begin{abstract}
This paper proposes novel approaches for designing control Lyapunov functions (CLFs) for constrained linear systems. We leverage recent configuration-constrained polyhedral computing techniques to devise piecewise affine convex CLFs. Additionally, we generalize these methods to uncertain systems with both additive and multiplicative disturbances. The proposed design methods are capable of approximating the infinite horizon cost function of both nominal and min-max optimal control problems by solving a single, one-stage, convex optimization problem. As such, these methods find practical applications in explicit controller design as well as in determining terminal regions and cost functions for nominal and min-max model predictive control (MPC). Numerical examples illustrate the effectiveness of this approach.
\end{abstract}

\end{frontmatter}

\section{Introduction}
\label{sec:introduction}
Control Lyapunov Functions (CLFs) have long been recognized as foundational tools in the analysis of stabilizability characteristics of linear and nonlinear systems, dating back to the pioneering works by Zubov~\cite{Zubov1965}, Artstein~\cite{Artstein1983}, and Primbs et al.~\cite{Primbs1999}. In particular, since their inception, piecewise affine CLFs and their polyhedral domains~\cite{Gutman1987} have received significant attention due to their flexibility and effectiveness in various control applications. By now, CLFs can be considered as a standard tool in control theory and a long list of methods, design strategies, and applications exist, as surveyed in~\cite{Giesl2015}.

In practical applications, however, the mere design of CLFs for stabilizing a system often falls short of the broader objectives. Instead, one is often tasked with the more intricate challenge of synthesizing CLFs that closely approximate the infinite horizon cost functions of optimal control problems. As discussed in~\cite{Bertsekas2012}, such approximately optimal CLFs are needed for achieving desired performance levels. Specifically, in the realm of model predictive control (MPC), CLFs are often used as terminal cost functions. Here, the objective is twofold: first, to leverage the inherent stability guarantees of CLFs; and second, to ensure that these functions accurately approximate the infinite horizon cost, thereby optimizing the closed-loop performance of the controller. This balance of computational tractability, stability, and optimality has been a focal point in the seminal work of Chen and Allg\"ower on quasi-infinite horizon MPC ~\cite{Chen1998}, as surveyed and generalized by Mayne et al.\cite{Mayne2000} and subsequent researchers, forming the basis for numerous MPC designs and implementations~\cite{Rawlings2009}.

Drawing parallels with the approximation of infinite horizon costs in traditional optimal control problems using CLFs, the design of min-max CLF functions serves a similar purpose in the realm of min-max optimal control problems. These robust CLFs find applications in approximating the infinite horizon cost associated with such problems. Existing techniques for constructing min-max CLFs often rely on the implementation of a min-max dynamic programming recursion, as pioneered by Witsenhausen~\cite{Witsenhausen1968}. Specifically, in the field of min-max model predictive control (MPC), researchers have developed a diverse array of iterative approximate dynamic programming approaches to design CLFs that cater to the specific requirements of this robust control framework. Works such as those by Boyd et.al.~\cite{Boyd1994}, Scokaert and Mayne~\cite{Scokaert1998}, Bemporad et.al.~\cite{Bemporad2003}, and Kerrigan and Maciejowski~\cite{Kerrigan2004} have laid the foundation for this research. This work has been further advanced in~\cite{Bjornberg2004} as well as in~\cite{Lincoln2006}, who developed approximate dynamic programming methods for designing piecewise-affine CLFs for min-max MPC.

Beyond these iterative and dynamic programming-based methods, convex and non-convex optimization-based polyhedral computing techniques have been proposed for constructing piecewise-affine CLFs. Notable examples include the work of Blanchini et al.~\cite{Blanchini1994,Blanchini1995,Blanchini1999a} and Raković~\cite{Rakovic2020,Rakovic2024}. A survey of related historical and recent approaches can be found in the overview articles~\cite{Giesl2015,Nguyen2018}.

Additionally, in recent years, learning-based frameworks have gained popularity, leading to the proposal of neural networks with Rectified Linear Units (ReLU) for designing piecewise affine and other types of CLF functions~\cite{Chen2022,Fabiani2022,Karg2020}. However, it is important to note that none of the aforementioned studies have proposed a direct CLF approximation of the solution to infinite horizon Hamilton-Jacobi-Bellman (HJB) equations or their min-max counterparts through a single-stage polytopic computing approach.

The challenge in designing piecewise affine CLFs via direct single-stage convex optimization lies in the necessity to rely on so-called double description methods, which involve enumerating the facets and vertices of the CLF's polyhedral epigraph. In this context, an impressive and early attempt was made by Jones and Morari in~\cite{Jones2010}, who proposed an implicit double description method for polyhedral CLF design. However, this early attempt relied on computationally expensive enumeration methods and did not result in a fully flexible convex optimization-based design. Nevertheless, as demonstrated in the current paper, it is feasible to modify and apply a recently proposed configuration-constrained polytopic computing technique~\cite{Houska2024,Villanueva2024} to develop such convex optimization-based CLF design techniques.

\subsection{Contribution}
The main contribution of this paper consists of novel polyhedral computing tools for constrained linear systems and min-max infinite horizon optimal control. These tools enable the design and optimization of parameterized polyhedral epigraphs of configuration-constrained piecewise-affine convex CLFs through the solution of single-stage convex optimization problems. In detail, we provide a computationally tractable and fully convex characterization of such parametric CLFs. The main theoretical results supporting this characterization are presented in Theorems~\ref{thm::epigraphHJB} and~\ref{thm::MixMaxHJB}. The main application of these CLFs is that they can be used to find approximate solutions to the HJB equations of both nominal as well as min-max infinite-horizon optimal control problems for linear systems with convex stage cost and constraints.

\subsection{Overview}
The structure of the paper is as follows.
\begin{itemize}
\item Section~\ref{sec::CLF} provides a concise review of the definition and key properties of CLFs.
\item Section~\ref{sec::Poly} proposes a novel class of piecewise affine CLFs with configuration-constrained polyhedral epigraphs.
\item Section~\ref{sec::HJB} explores the application of this new CLF parameterization for approximating solutions to HJB equations.
\item Section~\ref{sec::numeric} presents numerical case studies.
\item And, finally, Section~\ref{sec::conclusion} concludes the paper.
\end{itemize}

\subsection{Preliminaries and Notation}
Apart from standard notation, we use the symbol $\mathbb R_+^n$ to denote the non-negative orthant in $\mathbb R^n$. The symbol $\one \in \mathbb R_+^n$ is used to denote the vector whose coefficients are all equal to $1$. Moreover, throughout this paper, we adopt the basic notation and definitions from the field of polytopic computing~\cite{Fukuda2020}. Generally, a polyhedron is a (potentially empty) set of the form
\[
\left\{ \ x \in \mathbb R^n \ \middle| \ Fx \leq z \ \right\}.
\]
In this context, $F \in \mathbb R^{m \times n}$ is called the facet matrix and $z \in \mathbb R^m$ the facet parameter~\cite{Houska2024}. Moreover, bounded polyhedra are called polytopes, which are obtained for facet matrices $F$ for which $F x \leq 0$ implies $x=0$. Details about the basic properties of such polyhedra and polytopes can, for example, be found in~\cite{Ziegler1995}.

A function $M: \mathbb R^{n_x} \to \mathbb R \cup \{ \infty \}$ is called positive definite if $M(0) = 0$ and $M(x) > 0$ for all $x \neq 0$. We denote the epigraph of $M$ by
\[
\mathrm{epi}(M) \ \defeq \ \left\{ \ \left(
\begin{array}{c}
x \\
y
\end{array}
\right) \ \middle| \ M(x) \leq y \  \right\}.
\]
Throughout this paper, we call $M$ a proper compact convex function, if its domain $\operatorname{dom}(M)$ is non-empty and compact and $\mathrm{epi}(M)$ is a closed convex set. 

If $v_1,\ldots,v_m \in \mathbb R^n$ is a collection of points, 
\begin{equation*}
\mathrm{conv} \left( v_1,\ldots,v_m \right) \ \defeq \ \left\{ \sum_{i=1}^m \theta_i v_i \middle|
\begin{array}{l}
\theta_1, \ldots, \theta_m \geq 0 \\
\sum_{i=1}^m \theta_i = 1,
\end{array}
\right\}
\end{equation*}
denotes its convex hull. If $m$ is finite, $\mathrm{conv} \left( v_1,\ldots,v_m \right)$ is a polytope, whose vertices are a subset of these $m$ points. Finally, we introduce the symbol
\begin{equation*}
K_n \ \defeq \ \left\{ \ (0,\ldots,0,\theta)^\tr \ \middle| \ \theta \geq 0 \ \right\} \subseteq \mathbb R^n
\end{equation*}
to denote the proper convex cone that is generated by the last unit vector in $\mathbb R^n$. This notation is motivated by the fact that every proper compact convex function $M$ that is piecewise affine on a finite polytopic partition of its domain admits a vertex representation of its epigraph. That is, for $n = n_x + 1$, we have
\begin{equation}
\label{eq::Mvertex}
\operatorname{epi}(M) \ = \ \mathrm{conv}(v_1,\ldots,v_m) \oplus K_n\;,
\end{equation}
where $\oplus$ denotes the Minkowski sum. This follows from the Minkowski-Weyl theorem~\cite{Minkowski1989,Weyl1934}, recalling that $\operatorname{epi}(M)$ is a polyhedron whenever $M$ is a piecewise affine proper convex function. Throughout this paper, we refer to~\eqref{eq::Mvertex} as the vertex representation of such a function. Finally, the indicator function of a closed convex set $\mathbb X \subseteq \mathbb R^{n}$ is denoted by 
\begin{equation*}
    I_{\mathbb X}(x) := 
    \begin{cases}
    0 &\text{if }  x\in\mathbb{X} \\
    +\infty &\text{otherwise} \;.
    \end{cases} 
\end{equation*}

\section{Control Lyapunov Functions}
\label{sec::CLF}
This section reviews the basic definition of CLFs for both deterministic as well as uncertain constrained linear systems and highlights their importance in various applications.

\subsection{Constrained Linear Systems}
This paper is concerned with constrained linear systems with $n_x$ states and $n_u$ controls, 
\begin{equation}
\label{eq::system}
x_{k+1} = A x_k + B u_k \quad \text{subject to} \quad
\left\{
\begin{array}{l}
x_k \in \mathbb X \\
u_k \in \mathbb U \;.
\end{array}
\right.
\end{equation}  
Here, $\mathbb X \subseteq \mathbb R^{n_x}$ and $\mathbb U \subseteq \mathbb R^{n_u}$ denote, respectively, closed and convex state and control constraint sets. Moreover, we may assume---without loss of generality---that $B$ has full column rank, as redundant controls can be eliminated.

\begin{assumption}
\label{ass::blanket}
The sets $\mathbb X$ and $\mathbb U$ are closed and convex and such that $(0,0) \in \mathbb X \times \mathbb U$. Moreover, we assume that $\mathrm{rank}(B) = n_u \leq n_x$.
\end{assumption}

In many applications, the above system is equipped with a convex stage cost function $L$, as we are interested in the so-called control performance, $\sum_{k \in \mathbb N} L(x_k,u_k)$, typically over an infinite time horizon. We work with stage cost functions that satisfy the following assumption.

\begin{assumption}
\label{ass::Lnominal}
The function $L: \mathbb R^{n_x} \times \mathbb R^{n_u} \to \mathbb R$ is convex, non-negative, and satisfies $L(0,0) = 0$.
\end{assumption}

Note that Assumption~\ref{ass::Lnominal} is an assumption on the stage cost $L$ that is satisfied in many practical applications. In particular, it is satisfied when the control objective is the stabilization of a set-point or a more general set, with a quadratic cost function being a classical choice~\cite{Rawlings2009}. Assumption~\ref{ass::Lnominal} is, however, not satisfied for more general economic control objectives, where stabilization is not the main control objective. Nevertheless, in applications that satisfy a dissipativity condition and have a known storage function, it is possible to replace the economic stage cost function by an equivalent convex and positive semi-definite tracking cost~\cite{Faulwasser2018}. Additionally, we note that Assumption~\ref{ass::Lnominal} does not exclude the case that $L = 0$, which is a reasonable choice whenever no stage cost is given.

Finally, we recall the definition of control invariant sets~\cite{Blanchini2008}.
\begin{definition}
A compact convex set $X \subseteq \mathbb X$ with $0 \in X$ is called $\lambda$-control invariant~\footnote{Sometimes $\lambda$-control invariant sets with $0 \leq \lambda < 1$ are also called $\lambda$-contractive sets; see, e.g.,~\cite[Section~4.2.4]{Blanchini2008}.}, with $\lambda \in [0,1]$, if
\begin{equation*}
\forall x \in X, \ \exists u \in \mathbb U, \quad Ax+Bu \in \lambda X.
\end{equation*}
A $1$-control invariant set is called control invariant (CI).
\end{definition}

\subsection{Control Lyapunov Functions}
Control Lyapunov Functions (CLFs) constitute a fundamental tool for analyzing the stabilizability of control systems. In general settings, these functions are required to be merely lower semi-continuous and positive definite. However, in the context of linear-convex system theory, it is often sufficient to focus on convex CLFs. Consequently, in this article, we adhere to the following formal definition of (convex) $L$-CLFs.

\begin{definition}
\label{def::CLF}
A function $M: \mathbb R^{n_x} \to \mathbb R \cup \{ \infty \}$ is called an $L$-CLF if it is positive definite, proper compact convex, and satisfies the control Lyapunov inequality (CLI)
\footnote{\label{fn::Minimizer} If Assumptions~\ref{ass::blanket} and~\ref{ass::Lnominal} hold, the minimizer on the right side of~\eqref{eq::CLF} exists. This is because the map $$u \to L(x,u) + I_{\mathbb X}(x) + M(Ax+Bu)$$ is for every $x \in \mathbb R^{n_x}$ a compact, convex, and non-negative function. Here, the compactness statement follows from the requirement that $M$ is proper compact convex and our assumption that $B$ has full rank.}
\begin{equation}
\label{eq::CLF}
M(x)  \geq  \min_{u \in \mathbb U} \  L(x,u) + I_{\mathbb{X}}(x) + M(Ax+Bu) 
\end{equation}
for all $x\in\mathbb{R}^{n_x}$.
\end{definition}

\begin{remark}
\label{rem::StateConstraints}
In contrast to other existing definitions in the literature, state constraints of the form $Ax + Bu \in \mathbb X$ are not explicitly enforced in the CLI. They are, however, enforced implicitly: if $M$ satisfies the CLI, then we must have $M \geq I_{\mathbb X}$.
\end{remark}

Definition~\ref{def::CLF} contains the definition of convex control invariant sets as special case. Namely, a compact convex set is CI if and only if it is the sublevel set of a $0$-CLF, see~\cite{Blanchini2008}. Moreover, if the inequality in~\eqref{eq::CLF} is replaced with an equality, it coincides with the HJB equation of an infinite horizon optimal control problem with stage cost $L$, as discussed in~\cite{Bertsekas2012}. More precise statements regarding the relations between Definition~\ref{def::CLF}, infinite horizon optimal control, and HJBs can be found in Section~\ref{sec::HJB}.




\subsection{Robust Control Lyapunov Functions for Uncertain Linear Systems}
\label{subsec::rclf}
Definition~\ref{def::CLF} extends naturally to constrained uncertain linear systems of the form
\begin{equation}
\label{eq::system2}
x_{k+1} = A_k x_k + B_k u_k + w_k \quad \text{subject to} \quad
\left\{
\begin{aligned}
&x_k \in \mathbb X \\
&u_k \in \mathbb U,
\end{aligned}
\right.
\end{equation}
where the sequence $(A_k,B_k,w_k) \in \mathbb D$ is unknown. In this context, the uncertainty set,
\[
\mathbb D \subseteq \mathbb R^{n_x \times n_x} \times \mathbb R^{n_x \times n_u} \times \mathbb R^{n_x},
\]
is assumed to be non-empty, compact, and convex. In analogy with the nominal case, we recall the concept of robust control 
invariance.

\begin{definition}
\label{def::RCI}
A compact convex set $X \subseteq \mathbb X$ is called robust control invariant (RCI) if
\begin{equation*}
    \forall x \in X, \ \exists u \in \mathbb U, \forall (A,B,w) \in \mathbb D, \quad Ax+Bu + w \in X.
\end{equation*}
\end{definition}
In the following, we will additionally work with a non-standard extension of the above definition, which can be interpreted as a generalization of the concept of $\lambda$-contractiveness to uncertain control systems. The role of this definition in the ongoing developments will be explained in Section~\ref{sec::minmaxHJB}.
\begin{definition}
\label{def::lambdaRCI}
Let $X_{\sf s} \subseteq \mathbb X$ be a compact and convex RCI set and $\lambda \in [0,1)$ be a given constant. A compact convex set $X \subseteq \mathbb X$ is called robust $\lambda$-contractive relative to $X_{\sf s}
$ if
\begin{align*}
& \forall x \in X, \ \exists u \in \mathbb U, \forall (A,B,w) \in \mathbb D, \\
& \quad Ax+Bu + w \in \lambda X \oplus (1-\lambda) X_{\sf s},
\end{align*}
where $\oplus$ denotes the Minkowski sum
\footnote{The usual notion of robust $\lambda$-contractiveness corresponds to Definition~\ref{def::lambdaRCI} with $X_{\sf s} = \{0\}$, compare~\cite[Definition~4.19]{Blanchini2008}.}
.
\end{definition}

Uncertain linear systems are also often equipped with stage cost functions. While Assumption~\ref{ass::Lnominal} is sensible in the nominal case, slightly different assumptions
must be considered for uncertain systems. This is due to the fact that, in general, a system with additive uncertainties cannot be stabilized at a single point. Instead, the goal is to steer the system's state into a feasible RCI set $X_\mathrm{s} \subseteq \mathbb X$. 

\begin{definition}
\label{def::LPSD}
Let $X_{\mathrm{s}} \subseteq \mathbb X$ be a given RCI set. We say that a convex function $L: \mathbb R^{n_x} \times \mathbb R^{n_u} \to \mathbb R_+$ is positive semi-definite with respect to $X_{\mathrm{s}}$ for System~\eqref{eq::system2}, if
\begin{equation*}
\forall x\in X_{\sf s}, \ \exists u \in \mathbb U: \
\left\{
\begin{array}{l}
L(x,u) = 0 \quad \text{and} \\
\forall (A,B,w) \in \mathbb D,\\
Ax+Bu+w \in X_{\sf s}.
\end{array}
\right.
\end{equation*}
\end{definition}

The above definition is the basis for the following formal assumption.

\begin{assumption}
\label{ass::Lrobust}
The cost function $L: \mathbb R^{n_x} \times \mathbb R^{n_u} \to \mathbb R_{+}$ is convex and positive semi-definite with respect to at least one RCI set $X_\mathrm{s} \subseteq \mathbb X$.
\end{assumption}

Note that, as in the nominal case, Assumption~\ref{ass::Lrobust} does not exclude the case that we set $L = 0$. This assumption does, however, ensure that the following definition of min-max CLFs makes sense.

\begin{definition}
\label{def::MinMaxCLF}
A function $M: \mathbb R^{n_x} \to \mathbb R_+ \cup \{ \infty \}$ is called a robust $L$-CLF if it is non-negative, proper compact convex, satisfies $M(x_\mathrm{s})=0$ for at least one $x_\mathrm{s} \in \mathbb X$, and satisfies the robust CLI
\footnote{In contrast to the nominal case, a robust $L$-CLF function $M$ is not required to be positive definite. The existence of the minimum on the right side of the robust CLI is, however, still guaranteed as long as Assumptions~\ref{ass::blanket} and~\ref{ass::Lrobust} are satisfied. Similarly, the maximum 
exists since $\mathbb D$ is assumed to be compact.}
\begin{equation*}
M(x) \geq \min_{u \in \mathbb U} \max_{(A,B,w) \in \mathbb D} L(x,u) + I_{\mathbb{X}}(x) + M(Ax+Bu+w) 
\end{equation*}
for all $x \in \mathbb R^{n_x}$.
\end{definition}

In analogy to Definition~\ref{def::CLF},
Definition~\ref{def::MinMaxCLF} contains RCI sets as a special case. That is, a compact convex set is RCI if and only if it is the sublevel set of a robust $0$-CLF. Finally, the discrete-time min-max HJB for worst case optimal control on infinite time horizons can be obtained by replacing the inequality in the robust CLI with an equality; see~\cite{Bertsekas2012,Bjornberg2004,Witsenhausen1968}.

\section{CLF Construction via Polytopic Computing}
\label{sec::Poly}
This section introduces a necessary, sufficient, and computationally tractable
condition that characterizes the set of $L$-CLFs whose epigraph is a configuration-constrained polyhedron.

\subsection{Representation of Epigraphs}
Let $J: \mathbb R^{n_x} \to \mathbb R \cup \{ \infty \}$ be a proper compact convex function. Since its epigraph is a closed convex set, it can be approximated by a polyhedron sharing the same recession cone~\cite{Ney1995,Dorfler2022}. This is equivalent to approximating $J$ with a piecewise affine convex function $M$, whose domain is a polytope. One possible strategy is to first approximate $\mathrm{dom}(J)$ by a polytope of the form
\[
\left\{ \ x \in \mathbb R^{n_x} \ \middle| \  G_1 x \leq z_1 \ \right\},
\]
by selecting a facet matrix $G_1 \in \mathbb R^{\mathrm{f}_1 \times n_x}$ and the parameter $z_1$. The second step is to construct an augmented facet matrix 
\begin{equation}
\label{eq::F}
F \ = \ \left( G \ h \right) \ = \ 
\begin{pmatrix}
G_1 & 0 \\
G_2 & h_2
\end{pmatrix}
 \ \in \ \mathbb R^{\f \times n}
\end{equation}
with $n = n_x + 1$ and $\f = \mathrm{f}_1+\mathrm{f}_2$. Here, $G$ is obtained by augmenting the matrix $G_1$ with a matrix $G_2$ and a vector $h$ whose first $\mathrm{f}_1$ components are equal to zero, while the remaining $\mathrm{f}_2$ components are strictly negative.

\begin{assumption}
\label{ass::F}
The matrix $F\in\mathbb{R}^{{\sf f}\times n}$ satisfies~\eqref{eq::F} with a matrix $G_1\in\mathbb{R}^{\f_1\times n_x}$ such that $G_1 x \leq 0$ implies $x = 0$.\footnote{If the matrix $G_1$ is such that $G_1 x \leq 0$ implies $x = 0$, the set $\{ x \mid G_1 x \leq z_1 \}$ is bounded. This is necessary in practice, as finite piecewise affine approximations of convex functions on unbounded domains cannot be expected to be accurate.}
Moreover, 
$\f_2 > 0$ and the coefficients of $h_2\in\mathbb{R}^{\f_2}$ are strictly negative, that is $h_2 < 0$.
\end{assumption}

The matrix $F$ can be used to construct polyhedra 
\begin{equation*}
P(z) \ \defeq \ \left\{ \ 
\begin{pmatrix}
x \\
y
\end{pmatrix}\in\mathbb{R}^{n}
 \ \middle| \ Gx + h y \ \leq \ z \  \right\},
\end{equation*}
with parameter $z \in \mathbb R^{\f}$. Assumption~\ref{ass::F} ensures that such polyhedra can be used to represent epigraphs of proper compact convex functions.

\begin{proposition}
\label{prop::M}
Let $F$ satisfy Assumption~\ref{ass::F} and let $z$ be such that $P(z)$ is non-empty. Then, there exists a unique proper compact convex function $M: \mathbb R^{n_x} \to \mathbb R \cup \{ \infty \}$ such that $\mathrm{epi}(M) = P(z)$.
\end{proposition}

\begin{pf}
The set $P(z)$ is, by construction, a polyhedron and thus non-empty, convex, and closed. By Assumption~\ref{ass::F}, we have $h_2 < 0$. Thus, $P(z)$ must be the epigraph of a proper closed and convex function. This function is then also unique, as the given epigraph uniquely determines this function. Additionally, since $G_1 x \leq 0$ implies $x = 0$, the domain of this function is also bounded, implying the statement of the proposition.\hfill\hfill\qed
\end{pf}

\begin{remark}\label{rem::M}
If the conditions of Proposition~\ref{prop::M} are satisfied, the function $M$ can be constructed explicitly, by a direct comparison of $\operatorname{epi}(M)$ and $P(z)$. It is given by
\begin{equation*}
M(x) = \left\{ 
\begin{array}{cl}
\underset{i \in \{1,\ldots,{\sf f}_2\}}{\max} \ \dfrac{(z_2 - G_2 x)_{i}}{(h_2)_i} \ & \text{if} \ \ G_1 x \leq z_1 \\
\infty & \text{otherwise}.
\end{array}
\right.
\end{equation*}
\end{remark}

\subsection{Configuration Templates}
The class of polyhedra 
satisfying the conditions of Proposition~\ref{prop::M} is combinatorially large and also does not admit a jointly affine parameterization of facets and vertices~\cite{Houska2024}. As such, working with general polyhedra is not entirely practical for the development of efficient CLF design methods. Instead, we propose restricting the search to a specific set $\mathcal P(E)$ of polyhedra, which share the same face configuration. In detail, we propose the following definition of $\mathcal P(E)$, which is inspired by recent advancements in the field of polytopic computing~\cite{Villanueva2024}.

\begin{definition}
The triplet $(F,E,V)$ with facet matrix $F \in \mathbb R^{\f \times n}$, edge matrix $E \in \mathbb R^{\ep \times \f}$, and vertex matrices $V_1,V_2,\ldots, V_\vp \in \mathbb R^{n \times \f}$ is called a configuration triplet if
\begin{equation*}
\left.
\begin{array}{rcl}
P(z) &=& \{ x \in \mathbb R^{n+1} \mid Fx \leq z \} \ \text{and} \\
P(z) &=& \mathrm{conv} \left( V_1 z, \ldots, V_{\vp} z \right) \oplus K_n
\end{array}
\right\} \quad \Longleftrightarrow \quad Ez \leq 0,
\end{equation*}
with $z \in \mathbb R^\f$. The associated set of configuration-constrained polyhedra is denoted by
\begin{equation*}
\mathcal P(E) \ \defeq \ \left\{ \ P(z) \ \middle| \ \exists z \in \mathbb R^\f: \ E z \leq 0 \ \right\}.
\end{equation*}
\end{definition}

\begin{remark}\label{rem::triples}
Configuration triplets can be constructed in various ways. One option is, for example, to start with a non-empty template polytope 
$P(\overline z) = \{ \ x  \mid Fx \leq \overline z \ \}$, assumed to be simple and in minimal facet representation. In this case one can compute the vertex matrices $V_1,V_2,\ldots,V_\vp$ and the associated edge matrix $E$ as elaborated in~\cite[Section~3.5]{Villanueva2024}. If it happens that $P(\overline z)$ is not simple, one can add a small 
perturbation to the facet parameter such that $P(\overline z)$ is simple
after rescaling $F$, or directly use the considerations in~\cite[Proposition 1]{Houska2024a}. A survey of other efficient methods and tools for constructing such triplets can also be found in~\cite[Sections~2.10 and 2.13]{Houska2024}. There, it is also explained that $E$ can be constructed by enumerating the edges of the simple polyhedron $P(\overline z)$, motivating the name ``edge matrix" for $E$. Generally, it can be shown that the number $\ep$ of rows of $E$ is smaller than or equal to the number of edges of $P(\overline z)$; see~\cite[Section~2.13]{Houska2024}.
\end{remark}

If $(F,E,V)$ is a configuration triplet, if Assumption~\ref{ass::F} holds, and if $z$ satisfies $E z \leq 0$, then the conditions of Proposition~\ref{prop::M} are satisfied. This means that there exists a proper compact convex function $M$ for which $\mathrm{epi}(M) \ = \ P(z)$. Moreover, under the same conditions, the set of minimizers of
\begin{equation}
\label{eq::minpoint}
\min_{x,y} \ y \quad \  \text{s.t.} \quad \left(
\begin{array}{c}
x \\
y
\end{array}
\right) \in P(z)
\end{equation}
must be a face of $P(z)$~\cite{Ziegler1995}. Of course, in the most basic case in which~\eqref{eq::minpoint} has a unique minimizer, this set will simply be a vertex of $P(z)$. Thus, in this case, we can give this vertex the number $1$, recalling that the vertex enumeration of $P(z) \in \mathcal P(E)$ does not depend on the choice of $z$.

\begin{assumption}
\label{ass::vertex}
The triplet $(F,E,V)$ is a configuration triplet, such that the vertex $V_1 z$ is the unique minimizer of~\eqref{eq::minpoint} for any $z \in \mathbb R^{\f}$ with $E z \leq 0$.
\end{assumption}

Note that the above construction and assumptions are such that the following statement holds.

\begin{proposition}
\label{prop::P2}
Let Assumption~\ref{ass::F} and~\ref{ass::vertex} hold. Then the following statements are equivalent.
\begin{enumerate}
\addtolength{\itemsep}{4pt}
\item The polyhedron $P(z) \in \mathcal P(E)$ is the epigraph of a positive definite proper compact convex function $M$.

\item We have $E z \leq 0$ and $V_1 z = 0$.

\end{enumerate}
\end{proposition}

\begin{pf}
As mentioned above, the statement of this proposition is a direct consequence of Proposition~\ref{prop::M}. Here, the constraint $E z \leq 0$ is equivalent to the condition $P(z) \in \mathcal P(E)$. And, since $Ez \leq 0$ holds, the positive definiteness of the function $M$ is equivalent to the constraint $V_1 z = 0$, as guaranteed by Assumption~\ref{ass::vertex}. \qed
\end{pf}

\subsection{CLFs with Polyhedral Epigraphs}

The following theorem establishes computationally tractable convex feasibility conditions that characterize the set of all $L$-CLF functions with bounded domain, whose epigraph is a polyhedron with a given face configuration.

\begin{theorem}
\label{thm::epigraphHJB}
Let Assumptions~\ref{ass::blanket},~\ref{ass::Lnominal},~\ref{ass::F}, and~\ref{ass::vertex} be satisfied and let the matrices $R_i$ and $s_i$ be such that  $V_i = (R_i^\tr,s_i)^\tr$. Then, the following statements are equivalent.

\begin{enumerate}
\addtolength{\itemsep}{4pt}
\item The polyhedron $P(z) \in \mathcal P(E)$ is the epigraph of an $L$-CLF function.

\item We have $E z \leq 0$, $V_1 z = 0$, $R_i z \in \mathbb X$, and there exist control inputs $v_i \in \mathbb U$ and a vector $y \in \mathbb R^{\sf v}$ such that
\begin{align*}
L(R_i z,v_i) + y_i & \ \leq \  s_i^\tr z \\[0.16cm]
\text{and} \qquad  G A R_i z + G B v_i + h y_i & \ \leq \  z
\end{align*}
hold for all $i \in \{ 1,\ldots, \vp \}$.
\end{enumerate}
\end{theorem}

\begin{pf}
 First, Proposition~\ref{prop::P2} establishes the fact that $P(z) \in \mathcal P(E)$ is the epigraph of a positive definite and proper compact convex function $M$ if and only if $E z \leq 0$ and $V_1 z = 0$. Thus, it is left to show that the remaining conditions in the second statement are equivalent to enforcing the descent condition encoded in the CLI. The proof of this statement is divided in three parts. 

\texttt{Part I:} Since the right-hand expression in the CLI
\begin{equation*}
M(x) \geq \min_{u\in\mathbb{U}} \ L(x,u) + \mathbb{I}_{\mathbb X}(x) + M(Ax+Bu)
\end{equation*}
is convex, it is sufficient to check the descent condition at the points $R_i z$, which correspond to the projections of the vertices of $\operatorname{epi}(M)$ onto the domain of $M$. This implies that the CLI holds if and only if there exist control inputs $v_i \in \mathbb U$ with
\[
M(R_i z) \ \geq \ L(R_i z,v_i) + I_{\mathbb X}(R_i z) + M( A R_i z + B v_i )
\]
for all vertex indices $i \in \{ 1, \ldots, \vp \}$. Note that the derivation of the latter inequality can be regarded as a generalization of the key derivation step for the vertex control conditions for control invariant polytopes, which were originally introduced by Gutman and Cwickel~\cite{Gutman1986,Gutman1987}. Moreover, we use the equivalence
\[
\mathrm{dom}(M) \subseteq \mathbb X \qquad \Longleftrightarrow \qquad
\left\{
\begin{array}{l}
\forall i \in \{ 1, \ldots, \vp \}, \\
R_i z \in \mathbb X,
\end{array}
\right.
\]
to eliminate the indicator function from the inequality above by enforcing the constraints $R_i z \in \mathbb X$ explicitly.

\texttt{Part II:} Since we have $V_{i} = (R_i^\intercal, s_i)^\intercal$, the relation $M(R_i z) \ = \ s_i^\tr z$ holds for all $z$ with $E z \leq 0$. Thus, together with the consideration in the first part of the proof, we find that 
the CLI holds if and only if there exist vectors $v_i \in \mathbb U$ with
\[
s_i^\tr z \ \geq \ L(R_i z,v_i) + M( A R_i z + B v_i )
\]
and $R_i z \in \mathbb X$ for all vertex indices $i \in \{ 1, \ldots, n_\mathrm{v}\}$.

\texttt{Part III:} It remains to use the fact that the inequality $M(AR_i z + B v_i) \leq y_i$ holds if and only if
\[
(G,h) \left(
\begin{array}{c}
A R_i z + B v_i \\
y_i
\end{array}
\right) \ \leq \ z.
\]
Notice this is an immediate consequence of the fact that
$P(z) = \operatorname{epi}(M)$. Thus, by writing out this inequality and substituting $y_i$ in the inequality from the second part of the proof, the statement of the theorem follows. \hfill\hfill\qed
\end{pf}

\subsection{Robust CLFs for Min-Max Control}
\label{sec::minmaxCLF}
The developments from the previous sections can be extended for min-max control tasks for systems with given compact and convex uncertainty set~$\mathbb D$. To facilitate this extension, it is helpful to introduce the worst-case robust counterpart functions
\begin{equation*}
D_{i,j}(z,v) \ \defeq \ \max_{(A,B,w) \in \mathbb D} \ G_j A R_i z + G_j B v + G_j w
\end{equation*}
for all $i \in \{ 1, \ldots, \vp \}$ and $j \in \{ 1, \ldots, \f \}$. The functions $D_{i,j}$ are well-defined and convex, since we assume that $\mathbb D$ is compact and convex, recalling that the maximum over affine functions is convex. In the following, however, we simply assume that these functions are given, as they can be worked out explicitly for typical choices of $\mathbb D$ using known results from robust convex optimization~\cite{Bertsimas2011}.

\begin{example}
\label{ex::D}
Let us consider the practically relevant case that the uncertainty set $\mathbb{D} = \mathbb{W}_1 \times \mathbb{W}_2$ is a direct product of the polytopes
\begin{align*}
\mathbb W_1 &= \mathrm{conv}((A_1,B_1),\ldots,(A_{\bar l},B_{\bar l})) \\
\text{and} \quad \mathbb W_2 &= \{ w \in \mathbb R^{n_x} \mid Gw \leq \overline w \},
\end{align*}
where $\bar l \in \mathbb N$ denotes the number of vertices of the matrix polytope $\mathbb W_1$, while $\mathbb W_2$ is given in facet representation. Then, the robust counterpart functions are given by 
\begin{equation*}
D_{i,j}(z,v) = \max_{l \in \{ 1, \ldots, \bar l \}} \  G_j A_l R_i z + G_j B_l v +\overline{w}_j
\end{equation*}
for all $i \in \{ 1, \ldots, \vp \}$ and $j \in \{ 1, \ldots, \f \}$; see~\cite{Bertsimas2011,Houska2024}.
\end{example}

An additional consideration is that for uncertain systems it is often impossible
to find a robust $L$-CLF that is positive definite with respect to 0. Instead, 
one needs to work with configuration templates that explicitly permit the set 
of minimizers of~\eqref{eq::minpoint} to be a non-singleton set. Specifically, 
a---preferably small---RCI set. Note that this discussion is closely related to the corresponding construction in Definition~\ref{def::LPSD} and~Assumption~\ref{ass::Lrobust}. Consequently, Assumption~\ref{ass::vertex} is, in general, not adequate. An adequate assumption must, instead, reflect the fact that one facet of $P(z)$ needs to model the zero-level RCI set, as achieved by the following construction.

\begin{assumption}
\label{ass::flat}
The triplet $(F,E,V)$ is a configuration triplet. Moreover, the 
last row of the matrix $F\in\mathbb{R}^{\f \times n}$ is equal to
the negative $n$-th standard unit vector in $\mathbb{R}^{1\times \f}$. 
That is, $F_\f = -e_n^\intercal $.
\end{assumption}

After replacing Assumption~\ref{ass::vertex} with Assumption~\ref{ass::flat}, in analogy to Proposition~\ref{prop::P2}, the following statement holds.

\begin{proposition}
    Let Assumptions~\ref{ass::F} and~\ref{ass::flat} hold. Then, the following statements are equivalent.
    \begin{enumerate}
    \addtolength{\itemsep}{4pt}
    \item The polyhedron $P(z) \in \mathcal P(E)$ is the epigraph of a non-negative proper compact convex function $M$ with $M(x_\mathrm{s}) = 0$ for at least one $x_\mathrm{s} \in \operatorname{dom}(M)$.

    \item We have $E z \leq 0$ and $z_{\f} = 0$.   
    \end{enumerate}
    \end{proposition}

Based on this modified requirement on the choice of the triplet $(F,E,V)$, we are now able to characterize robust $L$-CLFs, whose epigraphs are configuration-constrained polyhedra.

\begin{theorem}
\label{thm::MixMaxHJB}
Let Assumptions~\ref{ass::blanket},~\ref{ass::Lrobust},~\ref{ass::F}, and~\ref{ass::flat} be satisfied, let $\mathbb D$ be compact and convex, and let the matrices $R_i$ and $s_i$ be such that  $V_i = (R_i^\tr,s_i)^\tr$. Then, the following statements are equivalent.

\begin{enumerate}
\addtolength{\itemsep}{4pt}
\item The polyhedron $P(z) \in \mathcal P(E)$ is the epigraph of a robust $L$-CLF.

\item We have $E z \leq 0$, $z_\f = 0$, $R_i z \in \mathbb X$, and there exist control inputs $v_i \in \mathbb U$ and a vector $y \in \mathbb R^{\sf v}$ such that
\begin{eqnarray}
L(R_i z,v_i) + y_i & \ \leq \ & s_i^\tr z \\[0.16cm]
\label{eq::Dij}
\text{and} \qquad  D_{i,j}(z,v_i) + h_j y_i & \ \leq \ & z_j
\end{eqnarray}
for all $i \in \{ 1,\ldots, \vp \}$ and all  $j \in \{ 1,\ldots, \f \}$.
\end{enumerate}
\end{theorem}

\begin{pf} 
The proof of this theorem is essentially analogous to the proof of Theorem~\ref{thm::epigraphHJB} after replacing the vertex constraint $V_1 z = 0$ with the facet condition $z_\f = 0$. Apart from this, one needs to modify the argument from the third part of the proof of Theorem~\ref{thm::epigraphHJB}, which is, however, straightforward as the inequality
\begin{equation*}
    D_{i,j}(z,v_i) + h_j y_i  \ \leq \ z_j
\end{equation*}
holds for all $i \in \{ 1,\ldots, \vp \}$ and all  $j \in \{ 1,\ldots, \f \}$,
if and only if 
\begin{equation*}
    M(AR_i z + B v_i + w) \ \leq \ y_i
\end{equation*}
holds for all $(A,B,w) \in \mathbb D$. This is
by definition of $D_{i,j}$, since they are constructed as the robust counterparts of the rows of the above inequality. \hfill\hfill\qed
\end{pf}

\section{CLF Design via Convex Optimization}
\label{sec::HJB}
The conditions specified in Theorems~\ref{thm::epigraphHJB} and~\ref{thm::MixMaxHJB} are convex, enabling the application of convex optimization methods for designing $L$-CLFs. In theory, one could formulate arbitrary, ideally convex, design objectives. 

In practical applications, however, one is often interested in identifying CLFs that closely approximate the solutions of the infinite horizon (min-max) Hamilton-Jacobi-Bellman (HJB) equation. Therefore, the current section primarily focuses on design objectives that lead to CLFs approximating such HJB solutions.

\subsection{Hamilton-Jacobi-Bellman Equation}

Let $X \subseteq \mathbb X$ be a compact convex $\lambda$-control invariant set with $0 \in X$. Moreover, let $\mathcal{M}(X)$ denote the set of functions $M:\mathbb{R}^{n_x}\to{R}$ such that $M$ is an $L$-CLF in the sense of Definition~\ref{def::CLF} and $\operatorname{dom}(M) = X$. 

Consider an arbitrary probability measure $\mu$ on $X$ with full support, $\mathrm{supp}(\mu) = X$. The goal of the following lemma is to establish conditions under which the 
minimizer $M^\star$ of 
\begin{equation}
\label{eq::Mstar}
\min_{M\in\mathcal{M}(X)} \ \int_X M \, \mathrm{d}\mu
\end{equation}
is well-defined and solves the HJB equation on $X$. That is, we have
\begin{equation}
    \label{eq::HJB}
    M^\star(x) \ = \ \min_{u \in \mathbb U} \ L(x,u) + I_X(x) + M^\star(Ax+Bu) 
\end{equation}
for all $x\in X$. Note that in terms of interpreting the above discrete-time HJB equation, the $\lambda$-control invariant set $X \subseteq \mathbb X$, not the set $\mathbb X$, should be regarded as our state constraint set, since we are interested in finding a CLF function, whose domain is equal to $X$. By definition, all functions $M \in \mathcal M(X)$ satisfy $M + I_X = M$. Consequently, the comments from Remark~\ref{rem::StateConstraints} also apply here---after replacing $\mathbb X$ with $X$. This explains why enforcing the state constraint $Ax+Bu \in X$ in~\eqref{eq::HJB} is then not necessary.

\begin{lemma}
\label{lem::hjb}
Let Assumptions~\ref{ass::blanket} and~\ref{ass::Lnominal} hold, let $X \subseteq \mathbb X$ be a compact convex $\lambda$-control invariant set with $0 \in X$ and $0 \leq \lambda < 1$, let $\mu$ be any probability measure on $X$ with full support and let $L$ be positive definite.
Then, the minimizer $M^\star \in \mathcal M(X)$ of~\eqref{eq::Mstar} exists and is unique. Moreover, $M^\star$ satisfies the HJB equation~\eqref{eq::HJB} for all $x \in X$.
\end{lemma}

Hamilton-Jacobi-Bellman equations are among the most extensively studied equations in optimal control theory, as evidenced in the literature~\cite{Bardi1997,Bellman1957,Bertsekas2012}. The above mention relation between the optimization problem~\eqref{eq::Mstar} and the HJB has, however, received much less attention. Hence, for the sake of maintaining the current article self-contained and given that the statement of the aforementioned lemma has---to the best of our knowledge---not been presented in this specific form, which is required in our context, we include a brief proof below.

\begin{pf}
Assumption~\ref{ass::Lnominal} ensures that $L$ is convex and hence Lipschitz continuous on $X$. Since we assume that $X$ is $\lambda$ control invariant with $0 \leq \lambda < 1$, this is already sufficient to conclude that $\mathcal M(X)$ is non-empty, because the system is exponentially stabilizable on $X$. Next, for any given CLF function $M \in \mathcal M(X)$, the associated finite horizon cost functions $J_N^M: \mathbb R^{n_x} \to \mathbb R_+ \cup \{ \infty \}$, given by
\begin{equation*}
\begin{alignedat}{2}
J_N^M(x_0)   =
& \min_{x,u} \ &&\sum_{k=0}^{N-1} L(x_k,u_k) + M(x_N) \\
& \ \text{s.t.} \ &&\left\{
\begin{aligned}
&\forall k \in \{ 0, \ldots, N-1 \} \\
&x_{k+1} = A x_k + B u_k \\
&x_k \in X, \ u_k \in \mathbb U,
\end{aligned} \right .
\end{alignedat}
\end{equation*}
are well-defined for all $x_0 \in \mathbb R^{n_x}$ and all horizon lengths $N \in \mathbb N$. 
This follows from an inductive argument using the positive definiteness 
of $M$ and the fact that a minimizer of $L(x,u)+M(Ax+Bu)$ exists 
in $\mathbb U$, see Footnote~\ref{fn::Minimizer}. 
Let us summarize the three most important properties of these cost functions
\begin{enumerate}
\addtolength{\itemsep}{3pt}

\item $J_N^M \in \mathcal M(X)$ for all $N \in \mathbb N$ and all $M \in \mathcal M(X)$;

\item $J_{N+1}^M \leq J_{N}^M$ for all $N \in \mathbb N$, and all $M \in \mathcal M(X)$; and

\item the infinite horizon cost
\begin{equation*}
J_\infty(x) = \lim_{N \to \infty} J_N^M(x)
\end{equation*}
exists and is finite for all $x \in X$. Moreover, this limit does not 
depend on the choice of $M \in \mathcal M(X)$.

\end{enumerate}
Note that the first two properties of the functions $J_N^M$ are immediate consequences of Bellman's principle of optimality. Moreover, the (pointwise) existence of the infinite horizon limit follows from the second statement by applying the monotone convergence theorem. Since $M$ is proper compact convex with $M(0)=0$, and hence continuous on its domain, for all $M \in \mathcal M(X)$, the infinite horizon limit $J_\infty$ does not depend on $M$. Additionally,  $J_\infty$ satisfies the HJB equation~\eqref{eq::HJB}~\cite{Bertsekas2012,Rawlings2009}. Together with our assumption that $L$ is positive definite, this implies that we also have $J_\infty \in \mathcal M(X)$.

Now, let us assume that we can find an $M \in \mathcal M(X)$ for which
\begin{equation}
\label{eq::BetterSolution}
\int_X M \, \mathrm{d}\mu \ < \ \int_X J_\infty \, \mathrm{d}\mu.
\end{equation}
In this case, it follows from the definition of $J_0^M = M$ and the monotonicity of the cost function sequence that we have
\begin{equation*}
\int_X M \, \mathrm{d}\mu \ = \ \int_X J_0 \, \mathrm{d}\mu \ \geq \ \int_X J_\infty \, \mathrm{d}\mu \ > \ \int_X M \, \mathrm{d}\mu.
\end{equation*}
But this is a contradiction. Thus, $J_\infty$ must be a minimizer of~\eqref{eq::Mstar}. The fact that $J_\infty$ is not only a minimizer, but actually the unique minimizer of~\eqref{eq::Mstar} follows by using a very similar argument: if there was a another minimizer $M \neq J_\infty$, we can find at least one probability measure $\mu$ for which~\eqref{eq::BetterSolution} holds, which then also yields a contradiction. Thus, the statement of the lemma holds.
\qed
\end{pf}

\begin{remark}
Note that solutions to HJBs are typically analyzed under weaker assumptions compared to those in Lemma~\ref{lem::hjb}. This is due to the fact that Lemma~\ref{lem::hjb} seeks to establish a condition under which the HJB equation~\eqref{eq::HJB} not only possesses a solution but also guarantees that this solution is an $L$-CLF that satisfies the requirements from Definition~\ref{def::CLF}. For instance, consider the case where $L = 0$, rendering $M = I_X$ a trivial solution to the HJB equation~\eqref{eq::HJB}. However, if $X \neq \{ 0 \}$, this function is not positive definite and, consequently, does not qualify as an $L$-CLF. Such scenarios are excluded in Lemma~\ref{lem::hjb} by imposing the sufficient condition that $L$ is positive definite. Alternatively, it would suffice to require that we have $\sum_{k \in \mathbb N} L(x_k,u_k) > 0$ along all admissible state and control trajectories with $x_0 \neq 0$.
\end{remark}

\subsection{Optimal CLFs with Polyhedral Epigraph}
 
This section concerns the construction of approximate solutions for~\eqref{eq::Mstar} by minimizing a suitable objective over the set of $L$-CLFs whose epigraph is a configuration-constrained polyhedron. Here, we recall that a computationally tractable representation of this set is given  
by the conditions of Theorem~\ref{thm::epigraphHJB}. In general, the idea is to solve optimization problems of the form
\begin{equation}\label{eq::cmin}
\begin{alignedat}{2}
&\min_{z,y,v} \ && \sigma(z) \\
&\ \text{s.t.} && \left\{
\begin{aligned}
&\forall i \in \{ 1, \ldots, \vp \}, \\
&L(R_i z,v_i) + y_i \leq s_i^\tr z \\
&G_1 A R_i z + G_1 B v_i \leq \lambda z_1 \\
&G_2 A R_i z + G_2 B v_i + h_2 y_i \leq  z_2 \\
&Ez \leq 0, \ V_1 z = 0, \ v_i \in \mathbb U, \ R_i z \in \mathbb X,
\end{aligned}
\right.
\end{alignedat}
\end{equation}
where $\sigma: \mathbb R^\f \to \mathbb R$ denotes an $L$-CLF design objective and $\lambda \in [0,1]$ a contraction parameter. If $\sigma$ is chosen to be a convex function, Problem~\eqref{eq::cmin}
is a convex optimization problem with ${\sf f}+{\sf v}(1+n_u)$ variables and ${\sf v}(1+{\sf f}+n_{\mathbb X} + n_{\mathbb U})+{\sf e} + n_x$ constraints. Here, $n_\mathbb X$ and $n_\mathbb U$ denote, respectively, the number of convex constraints needed to respresent the state and control constraints, $x\in\mathbb X$ and $u\in\mathbb U$.

Note that, a relatively simple and effective heuristic for choosing $\sigma$ is to set
\[
\sigma(z) \ = \ c^\tr z,
\]
where $c\in\mathbb{R}^{\f}$ is a strictly negative weight vector, with $c < 0$. Maximizing the facet parameters $z_i$, each with weight $-c_i > 0$, amounts to one way of maximizing some measure of the size of the epigraph $P(z) \in \mathcal P(E)$ of the $L$-CLF that is represented by these constraints.

There are, however, a couple of details that require a more careful discussion. First of all, if one wishes to interpret the solution $P(z)$ of~\eqref{eq::cmin} as the epigraph of an $L$-CLF $M$ that approximates the minimizer of~\eqref{eq::Mstar}, one needs to keep in mind that the statement of Lemma~\ref{lem::hjb} in general only holds on $\lambda$-control invariant sets $X$ with $0 \leq \lambda < 1$. In this context, the constraints
\begin{equation}
\label{eq::lambdaCIcond}
G_1 A R_i z + G_1 B v_i \leq \lambda z_1
\end{equation}
in~\eqref{eq::cmin} enforce the domain of the $L$-CLF $M$,
\[
\mathrm{dom}(M) \ = \ X \ \defeq \ \left\{ \ x \in \mathbb X \, \middle| \ G_1 x \leq z_1 \right\},
\]
to be a $\lambda$-control invariant polytope. Then, a possible variant 
of~\eqref{eq::cmin} can be obtained by first freezing $z_1$ for which~\eqref{eq::lambdaCIcond} is feasible and, ideally, for which $\operatorname{dom}(M)$ is a large set and then minimize over $v$ and $z_2$ only. In this case, the domain $X$ is frozen and it can then indeed be justified to minimize $\sigma(z) = c_2^\tr z$ for an arbitrary strictly negative weight, for example, $c_2 = -\one$. This is because, in this case, the goal is to approximate the minimizer of~\eqref{eq::Mstar}, which is known to be independent of the choice of the probability measure $\mu$.

Finally, a computationally more demanding but also significantly more systematic method for constructing a suitable CLF design objective $\sigma$ proceeds by pre-computing the optimal values
\begin{equation}\label{eq::zetaN}
\begin{alignedat}{2}
\zeta_i^N \ = \ 
& \max_{x,u} \ && G_i x_0 + \sum_{k=0}^{N-1} h_i L(x_k,u_k) + h_i \overline M(x_N) \\
& \ \text{s.t.} \ &&\left\{
\begin{aligned}
&\forall k \in \{ 0, \ldots, N-1 \}, \\
&x_{k+1} = A x_k + B u_k, \\
&x_k \in \mathbb X, \ u_k \in \mathbb U, x_N \in \mathbb X
\end{aligned} \right.
\end{alignedat}
\end{equation}
for all $i \in \{ 1, \ldots, \f \}$ and a large but finite time horizon $N$. Here, $0 \leq \overline M \leq M^\star$ denotes a suitable convex under-estimator of $M^\star$. For example, we can set $\overline{M} = 0$, or, if the stage cost $L$ is quadratic, we can choose $\overline{M}$ as the quadratic cost function of the corresponding unconstrained infinite-horizon Linear Quadratic Regulator (LQR), which also underestimates $M^\star$. Since we assume that $h_i < 0$ for all $i\in\{f_1+1,\ldots,f\}$--compare Assumption~\ref{ass::F}, the above maximization problems are convex optimization problems that can be solved efficiently, even if $N$ is relatively large. Moreover, under the rather mild additional assumption that $\mathbb X$ and $\mathbb U$ are bounded, the above maxima are bounded, too. The corresponding polyhedron $P(\zeta^N)$ is then, by construction, for all $N \in \mathbb N$, an outer approximation of the epigraph of the exact infinite horizon cost function $M^\star$ that solves the HJB under the assumption of Lemma~\ref{lem::hjb}. Since any feasible point $z$ of~\eqref{eq::cmin} yields an inner approximation, $P(z)$,
\[
P(\zeta^N) \ \supseteq \ \mathrm{epi}(M^\star) \ \supseteq \  P(z),
\]
one can use the function
\begin{equation*}
\sigma(z) \ = \ \left\| {\sf W}(z-\zeta^N) \right\|_\infty\;,
\end{equation*}
with an invertible weight matrix ${\sf W}\in\mathbb{R}^{n\times n}$, as design objective in order to directly minimize the (weighted) maximum distance between the upper and the lower bound on the exact infinite horizon cost $M^\star$ in the given directions $(G,h)$.

\subsection{Approximate Solutions to Min-Max HJBs}
\label{sec::minmaxHJB}
An extension of the considerations from the previous sections for uncertain constrained linear control systems is possible by replacing the conditions from Theorem~\ref{thm::epigraphHJB} by the corresponding conditions of Theorem~\ref{thm::MixMaxHJB}. 

In order to generalize the corresponding constructions from the previous section, however, we first need to choose a compact and convex RCI set $X_\mathrm{s} \subseteq \mathbb X$ and then choose a stage cost function $L$ that satisfies the positive semi-definiteness requirements from Assumption~\ref{ass::Lrobust} for this choice of $X_{\sf s}$. Moreover, we assume that $X$ is a compact and convex set that is robust $\lambda$-contractive set  relative to the same RCI set $X_{\sf s}$, with $\lambda \in [0,1)$; see~Definition~\ref{def::lambdaRCI}. Finally, in analogy to the nominal case, we use the symbol 
$\widehat{\mathcal M}(X)$ to denote the set of robust $L$-CLFs with domain $X$.

The following lemma establishes conditions under which the optimization problem
\begin{equation}\label{eq::robustMstar}
\min_{M \in \widehat{\mathcal M}(X)} \ \int_{X} M \, \mathrm{d}\mu
\end{equation}
admits a unique minimizer $M^\star$ that solves the min-max HJB equation on $X$,
\begin{align}
&M(x) = \notag \\
&\min_{u \in \mathbb U} \ \max_{(A,B,w) \in \mathbb D} \ L(x,u) + I_X(x) + M(Ax+Bu+w)
\label{eq::minmaxHJB}
\end{align}
for all $x\in X$. As for the nominal case, the set $X$ should here be regarded as a state constraint set. It needs to be taken into account in the HJB equation, as we are interested in designing CLFs with domain $X$.

\begin{lemma}
\label{lem::minmaxhjb}
Let $X_{\sf s} \subseteq \mathbb X$ be a given compact and convex RCI set, let Assumptions~\ref{ass::blanket} and~\ref{ass::Lrobust} hold for this choice of $X_{\sf s}$, and let $X \subseteq \mathbb X$ be a compact and convex set that is robust $\lambda$-contractive with respect to $X_{\sf s}$. Then the minimizer $M^\star\in\widehat{\mathcal M}(X)$ of~\eqref{eq::robustMstar} exists, is unique, and satisfies the min-max HJB equation~\eqref{eq::minmaxHJB} for all $x \in X$. This statement holds for any probability measure  $\mu$ on $X$ with full support.
\end{lemma}

\textit{Proof:} The proof of this lemma is over large parts analogous to the proof of Lemma~\ref{lem::hjb}. Nevertheless, one difference is that, for reasons that have been discussed in Section~\ref{subsec::rclf}, robust CLF functions are not required to be positive definite. As such, the positive semi-definiteness requirement on $L$ in Assumption~\ref{ass::Lrobust} is in this case already sufficient for our purposes. In detail, since $L$ is convex, it is also Lipschitz continuous. Moreover, as $X$ is robust $\lambda$-contractive, there exists for any $x_0 \in X$ a sequence of control laws $\kappa_k: \mathbb X \to \mathbb U$ that keeps the state $x_k$ of the uncertain closed-loop control system
\[
x_{k+1} = A x_k + B \kappa_k(x_k) + w_k
\]
in the sets $X_k = \lambda^k X \oplus (1-\lambda^k) X_{\sf s}$. For example, if $\kappa: \mathbb X \to \mathbb U$ denotes a control law that can be used to steer the system safely from $X$ to $\lambda X \oplus (1-\lambda)X_{\sf s}$ and if $\kappa_{\sf s}: \mathbb X \to \mathbb U$ can be used to keep the system's state within $X_{\sf s}$, then $\kappa_k = \lambda^k \kappa + (1-\lambda^k) \kappa_{\sf s}$ is such a sequence of control laws. This statement is a direct consequence of the condition from Definition~\ref{def::lambdaRCI}. This means that we have $x_k \in X_k$ for all $k \in \mathbb N$, independent of the realization of the uncertainty sequence. In other words, the distance of the system's state to the target set $X_{\sf s}$ is robustly exponentially stabilizable on $X$, which implies---together with Lipschitz continuity and positive semi-definiteness of $L$---that $\widehat{M}(X)$ is non-empty. In this context, another principal difference to the nominal case is that we define the cost functions $J_N^M$ by solving min-max (see e.g.~\cite{Scokaert1998,Kerrigan2004}) instead of nominal optimal control problems. 
In detail, we formulate the min-max optimal control problem
\begin{equation}
\label{eq::JNMrobust}
\begin{alignedat}{2}
    J_N^M(x_0)   =  
     \min_{x,u} \ &\max_{\ell \in \mathcal{I}_N} \ &&\sum_{k=0}^{N-1} L(x^{\ell}_k,u^{\ell}_k) + M(x^{\ell}_N) \\
    & \ \text{s.t.} \ &&\left\{
    \begin{aligned}
    &\forall k \in \{ 0, \ldots, N-1\}, \forall \ell \in \mathcal I_N, \\
    &x^{\ell}_{k+1} = A_{k}^\ell x^{\ell}_{k} + B^{\ell}_{k} u_{j(k,\ell)} + w^{\ell}_{k}, \\
    &x^{\ell}_k \in X, \ u_{j(k,\ell)} \in \mathbb U, \ x^{\ell}_{0} = x_0\;.
    \end{aligned} \right .
\end{alignedat}
\end{equation}
Here, $\ell$ indexes admissible uncertainty sequences, such that we have
\begin{equation*}
\bigcup_{\ell \in \mathcal I_N} ((A^{\ell}_0,B^{\ell}_0,w^{\ell}_{0}),\dots,(A^{\ell}_{N-1},B^{\ell}_{N-1},w^{\ell}_{N-1})) \ = \ \mathbb{D}^N,
\end{equation*}
where $\mathcal I_N$ denotes the (possibly uncountable) set of all such indices. Moreover, the control vectors are enumerated by the index $j(k,\ell)\in\mathbb{N}$ with $j(k,\ell) = j(k,i)$ whenever
\begin{equation*}
    \begin{aligned}
&((A^{\ell}_0,B^{\ell}_0,w^{\ell}_0),\ldots,(A^{\ell}_{k-1},B^{\ell}_{k-1},w^{\ell}_{k-1})) \\ 
&\quad =  ((A^{i}_0,B^{i}_0,w^{i}_0),\ldots,(A^{i}_{k-1},B^{i}_{k-1},w^{i}_{k-1}))
    \end{aligned}
\end{equation*} 
and $j(k,\ell) \neq j(k,i)$ otherwise, for all $i,\ell \in \mathcal I_N$ with $i\neq \ell$. The index notation in~\eqref{eq::JNMrobust} was introduced in~\cite{Scokaert1998} in the context of min-max model predictive control. In this context, the optimization variables, $x_k^\bullet$ and $u_{\bullet}$ are functions on $\mathcal I_N$, rendering~\eqref{eq::JNMrobust} a non-trivial functional optimization problem. Despite this complication, 
the existence of a min-max solution can be guaranteed under the listed assumptions. This can be proven using a standard min-max dynamic programming argument, see~\cite{Bertsekas2012}. Namely, by a simple induction over $N$, it follows that the functions $J_N^M$ are compact and convex. Consequently, as $B$ has full-rank, the mix-max dynamic programming recursion for $J_N^M$ can be used to establish the existence of min-max solutions. Thus, the 
finite-horizon cost functions $J^M_N$ are, as in the nominal case, 
well defined for all $x_0\in\mathbb{R}^{n_x}$ and all horizon 
lengths $N\in\mathbb{N}$.  

The argument of the proof of Lemma~\ref{lem::hjb} can now be recovered step 
by step as neither Bellman's principle of optimality nor the monotonicity argument used are affected by this change of definition, implying that the 
limit $J_\infty$ exists, does not depend on $M$, and satisfies the min-max 
HJB equation~\cite{Diehl2004}. The last contradiction step can be applied analogously, establishing the statement of the Lemma.\qed

\subsection{Optimal Robust CLFs with Polyhedral Epigraph}

In order to generalize the CLF design procedures from the previous sections to the min-max problem setting, we may assume that a configuration-constrained RCI polytope,
\begin{equation}\label{eq::RCIXs}
X_{\sf s} \ = \ \{ x \in \mathbb R^{n_x} \mid G_1 x \leq z^\mathrm{s} \} \subseteq \mathbb X,
\end{equation}
with $E z^\mathrm{s} \leq 0$, is given. Methods for computing such configuration-constrained polytopic target sets are surveyed in~\cite{Houska2024}. Next, an associated configuration-constrained polyhedral robust CLF approximation of the solution of the min-max HJB~\eqref{eq::minmaxHJB} can be found by solving the convex optimization problem
\begin{equation}
\label{eq::cminmax}
\begin{alignedat}{2}
&{\!\min_{z,y,v}} \ && \sigma(z) \\
&\text{s.t.}&& \left\{
\begin{aligned}
&\forall i \in \{ 1, \ldots, \vp \}, \\
&L(R_i z,v_i) + y_i \leq s_i^\tr z, \\
&Ez \leq 0, \ v_i \in \mathbb U, \ R_i z \in \mathbb X, \ z_\f = 0, \\
&\forall j \in \{ 1, \ldots, \f_1 \}, \\
&D_{i,j}(z,v_i) \leq \lambda z_j + (1-\lambda) z_j^\mathrm{s}, \\
&\forall j \in \{ \f_1 + 1, \ldots, \f_1 + \f_2 \}, \\
&D_{i,j}(z,v_i) + h_j y_i \leq  z_j.
\end{aligned}
\right.
\end{alignedat}
\end{equation}
The constraints $D_{i,j}(z,v_i) \leq \lambda z_j + (1-\lambda) z_j^\mathrm{s}$ enforce the $\lambda$-contractivity condition from Definition~\ref{def::lambdaRCI}, as elaborated in~\cite{Houska2024a}. Similar to the nominal case, the above problem has ${\sf f}+{\sf v}(1+n_u)$ optimization variables. In general, however, computing CLFs for min-max problems is more expensive than in the nominal case, as the number of convex inequality constraints in~\eqref{eq::cminmax} is equal to ${\sf v}(1+{\sf f }n_{\mathbb D}+n_\mathbb{U}+n_{\mathbb X}) +{\sf e}+{\sf f}$. Here, $n_{\mathbb D}$ denotes the number of constraints needed to evaluate the functions $D_{i,j}$. For example, if our uncertainty set $\mathbb D$ is chosen as in Example~\ref{ex::D}, we have $n_{\mathbb D} = \bar l$, recalling that $\bar l$ denotes the number of vertices of the matrix polytope.

The discussion regarding the design objective $\sigma$ parallels the one in the previous section. The primary difference lies when one considers minimizing the weighted maximum distance, $\sigma(z) = \| W(z-\zeta^N)\|_\infty$ between the epigraph $P(z)$ and an outer approximation $P(\zeta^N)$ of the epigraph of an optimal $L$-CLF that solved the min-max HJB equation. In such a case, the components of $\zeta^N$ are computed as 
\begin{equation}\label{eq::zetaNminmax}
\begin{alignedat}{2}
\zeta_i^N = \max_{x,u} & \min_{\ell \in \mathcal I_N} \ &&G_i x_0^\ell + \sum^{N-1}_{k=0} h_i L(x_k^\ell,u^\ell_k) + h_i \overline{M}(x_N^\ell) \\  
&\ \text{s.t.} && \left\{
\begin{aligned}
&\forall k\in\{0,\ldots,N-1\},  \forall \ell \in \mathcal I_N,\\
& x_{k+1}^\ell = A_k^\ell x_k^\ell + B_k^\ell u_{j(k,\ell)} + w_{k}^\ell \\
&x_{k}^\ell \in\mathbb{X}, \ u_{j(k,\ell)}\in\mathbb{U}, \
x_{N}^\ell\in\mathbb{X}\; , 
\end{aligned}\right.
\end{alignedat}
\end{equation}
where we use same notation as in the proof of Lemma~\ref{lem::minmaxhjb} and $0 \leq \overline M \leq M^\star = J_\infty^0$ is a convex under-estimator of the optimal infinite horizon cost. Unfortunately, this optimization problem is rarely ever tractable, as $\mathcal I_N$ is typically uncountable. An exception is the case where $\mathbb D$ is a polytope with given vertices, where it is sufficient to enumerate the extreme scenarios, as elaborated in~\cite{Kerrigan2004}. Nevertheless, even in this case solving~\eqref{eq::zetaNminmax} remains computationally intractable for larger $N$ due to an exponentially large number of uncertainty scenarios. Thus, in summary, while Problem~\eqref{eq::cminmax} is computationally tractable for a broad class of practical min-max optimal control problems and can serve as a viable means to obtain $L$-CLF approximations of the optimal min-max infinite horizon cost, achieving an exact computational verification of the precision of these approximations remains challenging.

\section{Numerical Illustration}
\label{sec::numeric}
Throughout this section, we present numerical results for one nominal and for one min-max test problems, which are used to illustrate the $L$-CLF design techniques in this paper.

\subsection{First Test Problem}
\label{sec::Problem1}
Our first test problem corresponds to a modified version of a standard benchmark problem for robust control design that has originally been proposed in~\cite{Goulart2006}. In the nominal setting, using the same notation as in~\eqref{eq::system}, the system matrices are given by
\begin{equation}
\label{eq::AB}
 A = \begin{pmatrix}
    1 & 1 \\
    0 & 1
\end{pmatrix}\quad \text{and}\quad
B = \begin{pmatrix}
    \frac{1}{2} \\
    1
    \end{pmatrix}\;.
 \end{equation}
The state and control constraint sets are given by
 \begin{equation*}
\mathbb{X} = [-1,2]^2 \quad\text{and}\quad \mathbb{U} = \left[-\frac{1}{2}, \frac{1}{2}\right].
\end{equation*}
Moreover, the stage cost function $L:\mathbb R^2 \times \mathbb R \to \mathbb R$ is given by
$L(x,u) = \Vert x \Vert^2_{\sf Q} + \| u \|_{\sf R}^2$ with weights
\begin{equation}
\label{eq::QR}
{\sf Q} = \begin{pmatrix}
1  & 0 \\
0  & \frac{1}{10}
\end{pmatrix} \quad \text{and} \quad {\sf R} = \frac{1}{10}\;.
\end{equation}
Note that this problem setting is used below to illustrate the proposed nominal $L$-CLF design methods.

\subsection{Second Test Problem}
\label{sec::Problem2}
Our second test problem considers an uncertain system of the form~\eqref{eq::system2} with $\mathbb D = \{ (A,B) \} \times \mathbb W$. Here, the matrices $A$ and $B$ are not uncertain, but given by~\eqref{eq::AB}. Moreover, the uncertainty set $\mathbb W$ of the additive process noise is set to
\begin{equation*}
\mathbb W = \left\{ \begin{pmatrix} 1 \\ 1 \end{pmatrix} \omega \ \middle| \  
\omega \in \left[ -\frac{1}{40},\frac{1}{40}\right] \right\}.
\end{equation*}
Finally, in order to construct a suitable stage cost function, we first compute a polyhedral approximation $X_{\sf s}$ of the minimal robust positive invariant set for the closed loop controlled system $x_{k+1} = (A+BK)x_k + w_k$, with
\[
K = -[0.895, 1.367] \quad \text{such that} \quad (A+BK)X_{\sf s} + \mathbb W \subseteq X_{\sf s}.
\]
A configuration-constrained RCI polytope $X_{\sf s}$ of the form~\eqref{eq::RCIXs} can be constructed by solving a single linear program, as explained in~\cite[Remark 5]{Villanueva2024}. A more thorough discussion on the construction
of configuration-constrained RCI polytopes can also be found in~\cite[Section~4.4]{Houska2024}.
The final set-regulation stage cost function 
$L:\mathbb R^2 \times \mathbb R \to \mathbb R$ is then defined by 
\begin{equation*}
L(x,u) = \min_{\xi\in X_{\sf s}} \ \Vert x - \xi \Vert^2_{\sf Q} + \| u - Kx \|_{\sf R}^2\;,
\end{equation*}
where ${\sf Q}$ and ${\sf R}$ are selected as in~\eqref{eq::QR}. Note that this setting will be used below to illustrate the proposed robust $L$-CLF design methods.

\subsection{Construction of Configuration Templates}
In order to set up the proposed CLF design optimization problems, we first need to choose a configuration triplet $(F,E,V)$. As discussed in Remark~\ref{rem::triples}, if a facet matrix $F$ is given, such that Assumption~\ref{ass::F} is satisfied, an associated triplet $(F,E,V)$ can be constructed by enumerating the edges and vertices of a selected non-empty polyhedron $P(\overline z) = \{ x \mid F x \leq \overline z \}$ and then apply the existing methods from~\cite{Houska2024,Villanueva2024} to compute $E$ and $V$. We discuss several strategies for constructing $F$ and $\overline z$, all of which are heuristic:

\begin{enumerate}

\item[S1)] The most basic strategy proceeds by selecting the first $\f_1$ rows of $F$ from the boundary of the $(n-1)$-hemisphere
$\left\{ c \in \mathbb{R}^{n} \ \middle| \ \Vert c \Vert_2 = 1, \ c_n \leq 0 \right\}$ and the remaining $\f_2$ rows from its interior, in such a way that Assumption~\ref{ass::F} holds. Now, one option is to set $\overline z = \one$.

\item[S2)] Another option is to select $F$ as above but then set $\overline z = \zeta^N$. Here, we compute $\zeta_1^N,\ldots,\zeta_{\sf f}^N$ by solving~\eqref{eq::zetaN} or~\eqref{eq::zetaNminmax}. This has the advantage that the outer approximation $P(\zeta^N) \supseteq \operatorname{epi}(M^\star)$ defines the configuration template.

\item[S3)] Note that strategies for selecting the first $\f_1$ rows of $F$ and strategies for pre-computing a (robust) $\lambda$-control invariant polytope $X$ are discussed in~\cite[Section~3.3 and~3.5]{Houska2024}). Based on this preparation, one possible strategy for selecting the remaining $\f_2$ rows of $F$ and $\overline z$ proceeds by evaluating the finite horizon cost functions $J_N^{\overline M}$ at selected sample points $\xi_1,\ldots,\xi_{\f_2} \in X$. One can then choose $F$ and $\overline z$ such that the polyhedron $P(\overline z)$ corresponds to the epigraph of a continuous piecewise affine interpolation of the points
\[
(\xi_{1},J^{\overline M}_{N}(\xi_{1})),\ldots,(\xi_{{\sf f}_2},J^{\overline M}_{N}(\xi_{{\sf f}_2})).
\]
Different methods for computing such a piecewise affine interpolation are discussed in~\cite[Section 6.5.5]{Boyd2004} as well as in~\cite[Proposition 3.2]{Bemporad2006}.
\end{enumerate}
Note that the first strategy, which sets $\overline z = \one$, results in a generic configuration template that may not always be adequate for approximating optimal value functions. In contrast, the second and third strategies guarantee that the templates can capture at least some approximation of $J_N^{\overline M} \approx M^\star$. However, it should be clear that all strategies are heuristics whose effectiveness is problem-dependent.

\subsection{Nominal L-CLF}
Our first numerical illustration concerns the nominal problem setting from Section~\ref{sec::Problem1}. We apply Strategy~S3) from the previous section to construct a facet matrix with $\f_1 = 8$ and $\f_2 = 37$ by interpolating sampling points on a pre-computed control invariant set, with $N = 5$ and with $\overline M$ being the optimal LQR cost function for the unconstrained problem. This leads to a configuration triplet $(F,E,V)$ with ${\sf f} = 45$, ${\sf v} = 88$, and ${\sf e} = 132$. 

\begin{figure}[h!]
    \centering
        \includegraphics[width=0.95\linewidth]{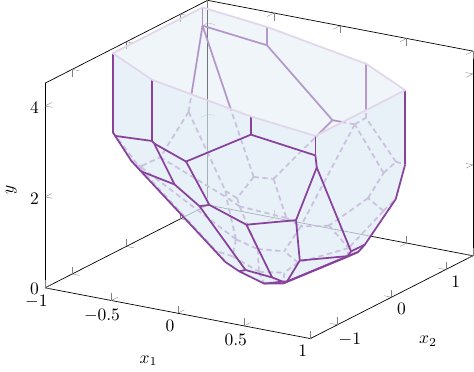}
        \caption{\label{fig::LCLF}Polyhedral 
    epigraph $P(z^\star)$ of the $L$-CLF function $M$.}
\end{figure}

Figure~\ref{fig::LCLF} shows $P(z^\star)$, the epigraph of the 
optimal $L$-CLF $M$ for the chosen template. Its facet
parameter $z^\star$ was computed by solving a convex quadratically 
constrained program with 
with $221$ decision variables and $4710$ constraints, based on Problem~\eqref{eq::cmin}, using the CLF design objective
\begin{equation*}
    \sigma(z) = \left\Vert {\sf W}\left(z - \zeta^{5}\right) \right\Vert_{\infty}\quad \text{with}\quad {\sf W} = 
    \begin{pmatrix}
    100 \cdot \mathbb{1}_8 & 0 \\
    0 & \mathbb{1}_{37}
    \end{pmatrix} 
\end{equation*}
with $\zeta^{5}$ computed as in~\eqref{eq::zetaN} with $\overline{M}$ defined as above. Here, $\mathbb{1}_{8}$ and $\mathbb{1}_{37}$ denote the unit matrices of 
sizes $8\times 8$ and $37\times 37$, respectively.
The optimization problem is sparse. For example, the vertex matrices have only between $7$ and $9$ (out of $135$) non-zero elements while the edge matrix has only $515$ (out of $5940$) non-zero elements. In fact, when 
aggregating all linear constraints, approximately $99.5\% $ of its entries are zeros. 

\begin{figure}[h!]
    \centering
        \includegraphics[width=0.95\linewidth]{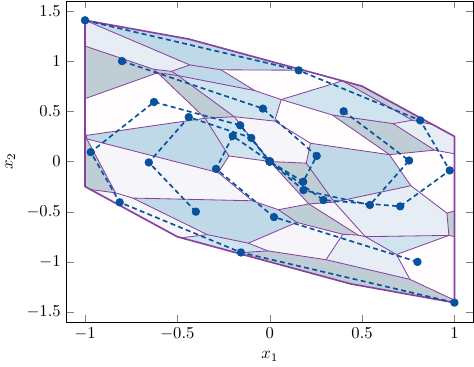}
        \caption{\label{fig::nominal-traj} Trajectories of the closed-loop system under $\mu_{M}$ and polyhedral partition of $\operatorname{dom}(M)$.}
\end{figure}

Figure~\ref{fig::nominal-traj} shows the polyhedral partition of $\operatorname{dom}(M)$
together with selected trajectories of the closed-loop system 
under the feedback law $\mu_{M}:\operatorname{dom}(M)\to\mathbb{U}$, given by
\begin{equation}\label{eq::illustration-mu}
\mu_{M}(x) \ \in \ \operatorname*{argmin}_{u\in\mathbb{U}} \ L(x,u) + M(Ax+Bu) \;.
\end{equation}
The feedback law $\mu_{M}$ can be evaluated
directly, without solving the above problem, since we have 
$\mu(R_{i}z^\star)= u^\star_i$ for all $i\in\{1,\ldots,\sf v\}$. Thus
$\mu(x)$ can be evaluated by solving a 
point location problem for $x$ followed by a linear 
interpolation of the optimal vertex 
control inputs $u^\star_i$ for the corresponding region. 

Notice that the construction above is of similar nature as explicit model predictive controllers~\cite{Herceg2013}. In particular,~\eqref{eq::illustration-mu} provides 
an approximation to the optimal, infinite horizon, feedback law whose complexity, unlike that of an exact explicit MPC controller, can be bounded
a priori and does not directly depend on the prediction horizon.

\begin{figure}[h!]
    \centering
        \includegraphics[width=0.95\linewidth]{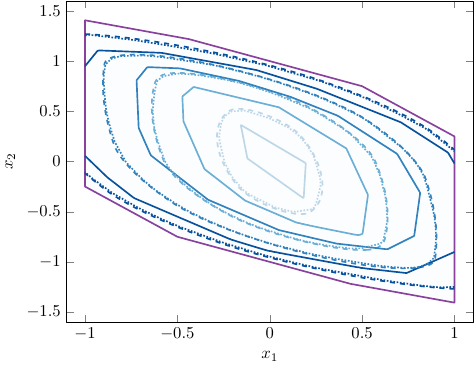}
        \caption{\label{fig::nominal-sublevel}(\legendbox{{189,215,231}} 0.1, \legendbox{{107,174,214}} 0.5, \legendbox{{49,130,189}} 1, and 
        \legendbox{{8,81,156}} 2)-Contours of $M$ (continuous lines), $M'$ (dotted lines), and $M^\star$ (dashed lines).}
\end{figure}

Finally, Figure~\ref{fig::nominal-sublevel} shows the $(0.1, 0.5,1,2)$-contours of the function $M$ (solid lines) and the function $M^\star$ solving the HJB equation (dashed lines). By construction, each sublevel set of $M$ is contained in the same sublevel set of 
$M^\star$ since $M(x) \geq M^\star(x)$ for all $x\in\operatorname{dom}(M)$. The maximum approximation error relative to the maximum value $M^\star_{\sf max}$ of $M^\star$ on its domain satisfies
\begin{equation*}
\frac{1}{M^{\star}_{\sf max}}\left(\max_{x\in\operatorname{dom(M)}}|M(x) - M^\star(x)|\right)\approx 0.41\;.
\end{equation*}
This error can be reduced, for example, by refining the polyhedral template. As
an illustration, a function $M'$ was constructed by using a different
configuration triple, which was generated by applying Strategy~S1) and randomly selected directions from the hemisphere. The dimensions of the corresponding random configuration triplet are ${\sf f} = 469$, ${\sf e} = 1401$, and ${\sf v} = 934$. Notice that the contour lines (dotted lines) of $M'$ approximately match those of $M^\star$. The maximum relative error of this approximation is approximately $0.08$.

\subsection{Robust L-CLF}
\label{sec::robustCLF}
This section illustrates the design of robust CLF functions based on the problem setting from Section~\ref{sec::Problem2}. For the sake of illustration, we apply Strategy S2) in this instance. In detail, this strategy is used to construct an initial facet matrix $F$, here with ${\sf f}_1 = 8$ and ${\sf f}_2 = 83$. Moreover, we have set $N = 7$ and $\overline M$ to be the LQR of the unconstrained LQR corresponding to the nominal system in order to solve~\eqref{eq::zetaNminmax}. Finally, this construction leads to a configuration triplet $(F,E,V)$ with ${\sf f}= 91$, ${\sf v}=180$, 
and ${\sf e}=270$.

\begin{figure}[h!]
    \centering
        \includegraphics[width=0.95\linewidth]{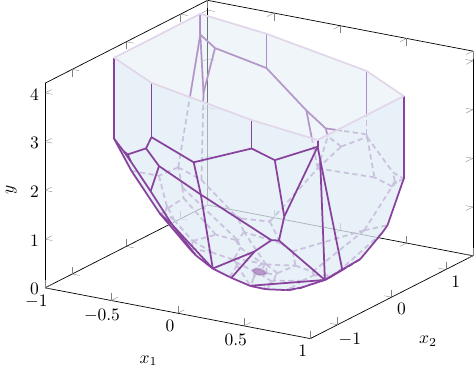}
        \caption{\label{fig::LCLF-robust}Polyhedral 
    epigraph $P(z^\star)$ of $M$, the robust $L$-CLF function. The dark (purple) facet satisfies $z^\star_{\sf f} = 0$.}
\end{figure}

Figure~\ref{fig::LCLF-robust} depicts the polyhedral epigraph $P(z^\star)$ 
of the optimal robust $L$-CLF $M$. The corresponding optimal facet parameter 
$z^\star$ was computed by solving the convex quadratically constrained program~\eqref{eq::cminmax}, with $451$ variables and $18002$ constraints.
As in the nominal setting, this problem is very sparse, too. For instance, the linear constraints
are encoded in a matrix with approximately $99.7\%$ zeros.

\begin{figure}[h!]
    \centering
        \includegraphics[width=0.95\linewidth]{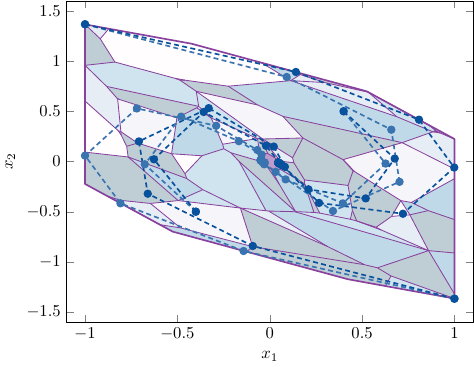}
        \caption{\label{fig::robust-traj} Trajectories of the closed-loop system under $\mu_{M}$ and polyhedral partition of $\operatorname{dom}(M)$. In (dark) purple, the set $\mathbb T$.}
\end{figure}

Figure~\ref{fig::robust-traj} shows the polyhedral partition of $\operatorname{dom}(M)$ and selected trajectories of the uncertain closed loop system under the feedback law
\begin{equation*}
\mu_{M}(x) \in \operatorname*{argmin}_{u\in\mathbb{U}} \ 
\max_{w\in\mathbb{W}} \ L(x,u) + M(Ax + Bu + w) \;.
\end{equation*}
As explained in the previous section, the feedback law can be evaluated based 
on the optimal vertex control inputs $v_{i}^\star$ from~\eqref{eq::cminmax}.
For each initial point, two trajectories are plotted corresponding to the extreme (but potentially not worst-case) 
uncertainty scenarios $w_{k} = \pm \frac{1}{40}$ for all $k \in \mathbb N$. 
Notice that, as predicted, all trajectories remain in $\operatorname{dom}(M)$ and converge to $\mathbb T$. As in the nominal case, the above constructions provides a means to construct an explicit controller with fixed complexity that approximates the 
optimal min-max infinite horizon optimal control law. 

\begin{figure}[h!]
    \centering
        \includegraphics[width=0.95\linewidth]{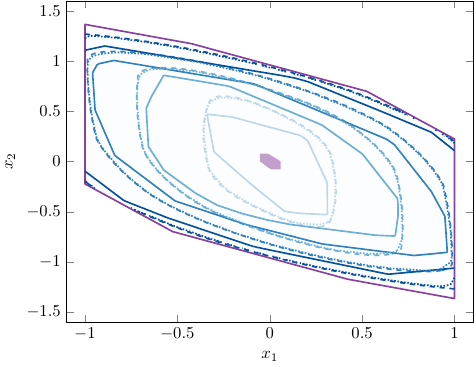}
        \caption{\label{fig::robust-sublevel}(\legendbox{{189,215,231}} 0.1, \legendbox{{107,174,214}} 0.5, \legendbox{{49,130,189}} 1, and 
        \legendbox{{8,81,156}} 2)-Contours of $M$ (continuous lines), $M'$ (dotted lines), and $M^\star$ (dashed lines). The 0-sublevel set (\legendbox{{195.5, 160, 206}}) corresponds, for all functions, to $\mathbb T$.}
\end{figure}

Finally, Figure~\ref{fig::robust-sublevel} shows selected contours of $M$, a very accurate numerical approximation of $M^\star$, and a second, finer, polyhedral 
approximation $M'$. This last  approximation was constructed using a configuration triple with ${\sf f} = 914$, ${\sf e} = 2487$, and ${\sf v} = 1658$, based on a random matrix $F$; see Strategy S1). Here, the relative approximation error of $M$ is approximately $0.47$, while that of $M'$ is approximately $0.11$.

\section{Conclusions}
\label{sec::conclusion}

This paper introduced a novel method for the design of control Lyapunov functions (CLFs) for constrained linear systems using configuration-constrained polyhedral computing techniques. The main contribution of this work is the development of convex conditions that characterize CLFs and robust CLFs with configuration-constrained polyhedral epigraphs, as detailed in Theorems~\ref{thm::epigraphHJB} and~\ref{thm::MixMaxHJB}. These conditions enable the construction of piecewise affine convex CLFs that approximate the infinite horizon value function for both nominal and min-max optimal control problems through the solution of a single convex optimization problem.

A notable feature of the proposed design techniques is that optimizing the epigraph of the polyhedral CLF inherently determines the associated control inputs at the vertices of the polyhedral domain partition. The resulting control inputs provide a computationally efficient approximation of the infinite horizon optimal control law, which can be evaluated using point location and convex interpolation techniques, as discussed and illustrated in Section~\ref{sec::numeric}.

The versatility of this approach holds promise for applications in determining terminal regions and cost functions in MPC. Moreover, it provides a valuable bridge between implicit and explicit control synthesis approaches. Specifically, the proposed methodology offers a means to construct an explicit controller that approximates the optimal infinite horizon feedback control law, with a complexity that is fixed a-priori. This characteristic ensures that the resulting controller is both computationally predictable and efficient, making it particularly well-suited for real-time applications, where fixed computational requirements are critical. The numerical examples presented in this paper highlight the practicality and effectiveness of the proposed CLF design strategies, underscoring their potential for wide-ranging applications in control theory and practice.

\bibliographystyle{plain}
\bibliography{references}

\end{document}